\def\ekv#1#2{\begeq\label{#1}#2\endeq}

\def\ably{arbitrarily}
\def\an{analytic}

\def\bdd{bounded}
\def\bdy{boundary}

\def\dop{differential operator}

\def\ev{eigenvalue}
\def\e{equation}
\def\fu{function}
\def\fy{family}

\def\hf{Hamilton field}
\def\hg{homogeneous}

\def\indep{independent}
\def\lhs{left hand side}
\def\mfld{manifold}

\def\neigh{neighborhood}
\def\nondeg{non-degenerate}
\def\op{operator}

\def\pb{problem}

\def\pol{polynomial}

\def\rhs{right hand side}
\def\sa{selfadjoint}

\def\sufly{sufficiently}
\def\sy{satisfy}

\def\Th{Theorem}

\def\ufly{uniformly}
\def\vf{vector field}
\def\wrt{with respect to}

\def\Im{{\rm Im\,}}
\def\diverg{\mathop{\fam0 div}}

\documentclass[10pt]{article}
\usepackage{amssymb,latexsym}
\title{The Calder{\'o}n problem with partial
data}
\author{Carlos E. Kenig\\Department of Mathematics\\
University of Chicago
Chicago, IL 60637, USA\\
and\\
Institute for Advanced Study\\
Princeton, New Jersey 08540, USA\\cek@math.uchicago.edu\\
\and
Johannes Sj{\"o}strand\\Centre de Math{\'e}matiques\\Ecole 
Polytechnique\\FR 
91120 Palaiseau, France\\johannes@math.polytechnique.fr\\ (UMR 7640,
CNRS)\and Gunther
Uhlmann
\\Department of Mathematics\\
University of Washington\\
Box 354350
Seattle, WA 98195, USA\\gunther@math.washington.edu}
\date{}
\def\wrtext#1{\relax\ifmmode{\leavevmode\hbox{#1}}\else{#1}\fi}

\def\begeq{\begin{equation}}
\def\endeq{\end{equation}}

\def\part#1{\frac{\partial}{\partial #1}}

\newcommand{\dist}{\mbox{\rm dist\,}}

\renewcommand{\Im}{\mbox{\rm Im\,}}

\newcommand{\supp}{\mbox{\rm supp}}

\newtheorem{dref}{Definition}[section]
\newtheorem{lemma}[dref]{Lemma}
\newtheorem{theo}[dref]{Theorem}
\newtheorem{prop}[dref]{Proposition}

\newtheorem{coro}[dref]{Corollary}

\newenvironment{proof}{\vspace{.3cm}\noindent{{\em 
Proof:}}}{\hfill$\Box$\medskip}
\begin{document}
\maketitle

\begin{abstract}
In this paper we improve an earlier result by Bukhgeim and Uhlmann 
\cite{BuUh}, by 
showing that in dimension $n\ge 3$, the knowledge of the Cauchy data for 
the Schr{\"o}dinger equation measured on possibly very small subsets of 
the 
\bdy{} determines uniquely the potential. We follow the general strategy 
of \cite{BuUh} but use a richer set of solutions to the Dirichlet \pb{}.
This implies a similar result for the problem of Electrical Impedance
Tomography which consists in determining the conductivity of a body
by making voltage and current measurements at the boundary.
\end{abstract}

\vskip 2mm
\noindent
{\bf Keywords and Phrases:} Dirichlet to Neumann map, Carleman estimates, 
analytic microlocal analysis. 

\vskip 1mm
\noindent
{\bf Mathematics Subject Classification 2000}: 35R30

\section{Introduction}\label{section0}
\setcounter{equation}{0}

 The Electrical Impedance Tomography (EIT) inverse problem consists in
determining the electrical conductivity of a body by making voltage and
current measurements at the boundary of the body. Substantial progress has
been made on this problem since Calder{\'o}n's pioneer contribution 
\cite{C}, and is also known as Calder{\'o}n's problem, in
the case where the measurements are made on the whole boundary. 
This
problem can be reduced to studying the Dirichlet-to-Neumann (DN) map
associated to the Schr{\"o}dinger equation. A key ingredient in several of
the results is the construction of complex geometrical optics for the
Schr{\"o}dinger equation (see \cite{U} for a survey).  Approximate 
complex geometrical optics 
solutions for the Schr\"odinger equation concentrated near planes are 
constructed in \cite{GrU} and concentrated near spheres in \cite{IsU}.

Much less is known if the DN map is only measured on part of the boundary. 
The only previous result that we are aware of, without assuming any 
a-priori condition on the potential besides being bounded, is in
\cite{BuUh}. 
It is shown there that if we measure the DN map restricted to, roughly 
speaking, a slightly more than half of the boundary 
then one
can determine uniquely the potential. The proof relies on a Carleman
estimate with an exponential weight with a linear phase. The Carleman
estimate can also be used to construct complex geometrical optics
solutions for the Schr{\"o}dinger equation. We are able in this paper to
improve significantly on this result. We show that measuring the DN map on
an arbitrary open subset of the boundary we can determine uniquely the
potential. We do this by proving a more general Carleman estimate
(Proposition 3.2) with exponential non-linear weights. This Carleman
estimate allows also to construct a much wider class of complex
geometrical optics than previously known (section 4).
We now state more precisely the main results.

In the following, we let $\Omega \subset\subset{\bf R}^n$, be an open 
connected set 
with $C^\infty $ \bdy{}. For the main results, we will also assume that 
$n\ge 3$. If $q\in L^\infty (\Omega )$, then we consider the \op{} 
$-\Delta +q:\, L^2(\Omega )\to L^2(\Omega )$ with domain $H^2(\Omega )\cap 
H_0^1(\Omega )$ as a \bdd{} perturbation of minus the usual Dirichlet 
Laplacian. $-\Delta +q$ then has a discrete spectrum, and we assume 
\ekv{0.1}{0\hbox{ is not an \ev{} of }-\Delta +q:\, H^2(\Omega )\cap H_0^1
(\Omega 
)\to L^2(\Omega ).
}
Under this assumption, we have a well-defined Dirichlet to Neumann map
\ekv{0.2}{{\cal N}_q:\, H^{1\over 2}(\partial \Omega )\ni v\mapsto 
{\partial _\nu 
u_\vert }_{\partial \Omega }\in H^{-{1\over 2}}(\partial \Omega ),
}
where $\nu $ denotes the exterior unit normal and $u$ is the unique 
solution in 
\ekv{0.3}
{H_\Delta (\Omega ):=\{ u\in H^1(\Omega );\,\, \Delta u\in L^2(\Omega )\} }
of the \pb{}
\ekv{0.4}
{
(-\Delta +q)u=0\hbox{ in }\Omega ,\ {u_\vert}_{\partial \Omega }
=v.  }
See \cite{BuUh} for more details, here we have slightly modified the 
choice 
of the Sobolev indices.

Let $x_0\in {\bf R}^n\setminus\overline{{\rm ch\,}(\Omega )}$, where ${\rm 
ch\,}(\Omega )$ denotes the convex hull of $\Omega $. Define the front and 
the back faces of $\partial \Omega $ by 
\ekv{0.5}
{
F(x_0)=\{ x\in\partial \Omega ;\, (x-x_0)\cdot \nu(x)\le 0\} ,\ B(x_0)=\{ 
x\in\partial \Omega ;\, (x-x_0)\cdot \nu(x) > 0\} .}
The main result of this work is the following:

\begin{theo}\label{Th0.1}
With $\Omega $, $x_0$, $F(x_0)$, $B(x_0)$ defined as above, let $q_1,q_2\in
L^\infty (\Omega )$ be two potentials satisfying (\ref{0.1}) and assume
that there exist open \neigh{}s $\widetilde{F},\,\widetilde{B}\subset \partial \Omega $ of
$F(x_0)$ and $B(x_0)\cup\{ x\in\partial \Omega ; (x-x_0)\cdot \nu =0\}$ respectively, such that
\ekv{0.6}
{
{\cal N}_{q_1}u={\cal N}_{q_2}u\hbox{ in }\widetilde{F},\hbox{ for all
}u\in H^{1\over 2}(\partial \Omega )\cap{\cal E}'(\widetilde{B}).} 
Then $q_1=q_2$.
\end{theo}

Notice that by Green's formula ${\cal N}_q^*={\cal N}_{\overline{q}}$. It 
follows that $\widetilde{F}$ and $\widetilde{B}$ can be 
permuted in (\ref{0.6}) and we get the same conclusion.

\par If $\widetilde{B}=\partial \Omega $ then we obtain the following result.

\begin{theo}\label{Th0.2}
With $\Omega $, $x_0$, $F(x_0)$, $B(x_0)$ defined as above, let 
$q_1,q_2\in 
L^\infty (\Omega )$ be two potentials satisfying (\ref{0.1}) and assume 
that there exists a \neigh{} $\widetilde{F}\subset \partial \Omega $ of 
$F(x_0)$, such that 
\ekv{0.7}
{
{\cal N}_{q_1}u={\cal N}_{q_2}u\hbox{ in }\widetilde{F},\hbox{ for all 
}u\in H^{1\over 2}(\partial \Omega ).
}
Then $q_1=q_2$.
\end{theo}

\par We have the following easy corollary,

\begin{coro}\label{Cor0.2} 
With $\Omega $ as above, let $x_1\in\partial \Omega $ be a point such that 
the tangent plane $H$ of $\partial \Omega $ at $x_1$ satisfies $\partial 
\Omega \cap H=\{ x_1\}$. Assume in addition, that $\Omega $
is strongly starshaped \wrt{} $x_1$. Let $q_1,q_2\in L^\infty (\Omega )$ 
and assume that there exists a \neigh{} $\widetilde{F}\subset \partial 
\Omega $ of $x_1$, such that (\ref{0.7}) holds. Then $q_1=q_2$.
\end{coro}

Here we say that $\Omega $
is strongly star shaped \wrt{} $x_1$ if every line through $x_1$ which 
is not contained in the tangent plane $H$ cuts the \bdy{} $\partial 
\Omega $ at precisely two distinct points, $x_1$ and $x_2$, and the 
intersection at $x_2$ is transversal.

Theorem \ref{Th0.1} has an immediate consequence for the Calder{\'o}n problem. 

Let $\gamma\in C^2(\overline \Omega)$
be a strictly positive function on $\overline \Omega$.
Given a voltage potential $f$ on the boundary,
the equation for the potential in the interior,
under the assumption of no sinks or sources of current in
$\Omega$,
is
$$
\diverg(\gamma \nabla u)=0 \hbox{ in } \Omega, \quad {u_|}_{\partial\Omega}=f.
$$
The Dirichlet-to-Neumann map is defined in this case as follows:
$$
{\cal N}_\gamma (f)= {\left(\gamma \partial_\nu  u\right)_
|}_{\partial\Omega}.
$$
It extends to a~bounded map
$$
{\cal N}_\gamma: H^{\frac{1}{2}}(\partial\Omega)\longrightarrow
H^{-\frac{1}{2}}(\partial\Omega).
$$

As a direct consequence of Theorem~1.1 we have 

\begin{coro}\label{coro0.3}
Let
$\gamma_i \in C^2(\overline\Omega)$, $i=1,2$,
be strictly positive. 
Assume that $\gamma_1=\gamma_2$ on $\partial\Omega$ and
$$
{\cal N}_{\gamma_1}u={\cal N}_{\gamma_2}u\mbox{ in }\widetilde{F},
\mbox{ for all }
u\in H^{1\over 2}(\partial \Omega )\cap{\cal E}'(\widetilde{B}).
$$
Then $\gamma_1=\gamma_2$.
\end{coro}

Here $\widetilde{F}$ and $\widetilde{B}$ are as in Theorem \ref{Th0.1}. 
It is well known (see for instance \cite{U}) that one can relate 
${\cal N}_\gamma$ and ${\cal N}_q$
in the case that $q=\frac{\Delta \sqrt\gamma}{\sqrt\gamma}$ with $\gamma >0$ 
by the formula
\begin{equation}\label{formula}
{\cal N}_q(f)= {(\gamma^{-\frac{1}{2}})_|}_{\partial \Omega} 
{\cal N}_\gamma(\gamma^{-\frac{1}{2}} f)+
\frac{1}{2} {\left(\gamma^{-1} \partial_\nu\gamma\right)_|}_{\partial\Omega} f.
\end{equation}
The Kohn-Vogelius result \cite{KV} implies that $\gamma_1=\gamma_2$  and
$\partial_\nu \gamma_1= \partial_\nu \gamma_2$
on $\widetilde{F}\cap\widetilde{B}$.
Then using (\ref{formula}) and Theorem \ref{Th0.1} we immediately get
Corollary \ref{coro0.3}.

\par
A brief outline of the paper is as follows. In section 2 we review the 
construction of weights that can be used in proving Carleman
estimates. In section 3
we derive the Carleman estimate (Proposition 3.2) that we shall use in the
construction of complex geometrical optics solutions for the 
Schr{\"o}dinger
equation. In sections 4, 
5 we use the Carleman estimate for solutions of
the inhomogeneous Schr{\"o}dinger equation vanishing on the boundary. This
leads to show that, under the conditions of Theorems 1.1 and 1.2, the 
difference
of the potentials is orthogonal in $L^2$ to a family of oscillating 
functions which
are real-analytic. For simplicity we first prove Theorem 1.2. In section 6 
we end the proof of Theorem 1.2 by choosing this family appropriately 
and using the wave
front set version of Holmgren's uniqueness theorem. Finally in section 7 we
prove the more general result Theorem 1.1.

\medskip 

\noindent{\it Acknowledgments.} The first author was supported in part by 
NSF
and at IAS by The von Neumann Fund, The Weyl Fund, The Oswald Veblen Fund 
and
the Bell Companies Fellowship. 
The second author was partly supported by the MSRI in Berkeley and the last
author was partly supported by NSF and a John Simon Guggenheim fellowship.

\section{Remarks about Carleman weights in the variable 
coefficient case}\label{SectionVC}
\setcounter{equation}{0}

\par In this section we review the construction of weights that can be 
used in proving Carleman estimates. The discussion is a little more 
general than what will actually be needed, but much of the section can be 
skipped at the first reading and we will indicate where.

Let $\widetilde{\Omega }\subset {\bf R}^n$, $n\ge 2$ be an open set, and 
let $G(x)=(g^{ij}(x))$ a positive definite real symmetric $n\times 
n$-matrix, depending smoothly on $x\in\widetilde{\Omega }$. Put 
\ekv{vc.1}
{p(x,\xi )=\langle G(x)\xi \vert \xi \rangle .}
Let $\varphi \in C^\infty (\widetilde{\Omega };{\bf R})$ with $\varphi '(x)\ne 
0$ everywhere, and consider
\ekv{vc.2}
{p(x,\xi +i\varphi '_x(x))=a(x,\xi )+ib(x,\xi ),}
so that with the usual automatic summation convention:
\ekv{vc.3}
{a(x,\xi )=g^{ij}(x)\xi _i\xi _j-g^{ij}(x)\varphi '_{x_i}\varphi '_{x_j},}
\ekv{vc.4}
{b(x,\xi )=2\langle G(x)\varphi '(x)\vert \xi \rangle =2g^{\mu \nu }\varphi 
'_{x_\mu }\xi _\nu .}
Readers, who are not interested in routine calculations, may go directly 
to the conclusion of this section.

\par 
A direct computation gives the Hamilton field $H_a=a'_\xi \cdot \partial_x
-a'_x \cdot \partial_\xi$ of
$a$:
\begin{eqnarray}\label{vc.4.5}
H_{a}=2g^{ij}(x)\xi _j\partial
_{x_i}-\partial _{x_\nu }(g^{ij})\xi _i\xi 
_j\partial _{\xi _\nu
}+\partial _{x_\nu }(g^{ij})\varphi '_{x_i}\varphi 
'_{x_j}\partial _{\xi _\nu
}+2g^{ij}\varphi ''_{x_i,x_\nu }\varphi 
'_{x_j}\partial _{\xi _\nu
},
\end{eqnarray} 
and 
\begin{eqnarray}\label{vc.5}
{1\over 2}H_{a}b&=&
2g^{ij}\xi _jg^{\mu \nu }\varphi ''_{x_i,x_\mu }\xi _\nu +2g^{ij}\xi_j
\partial _{x_i}(g^{\mu \nu })\varphi '_{x_\mu }\xi _\nu -\partial _{x_\nu
}(g^{ij})\xi _i\xi _jg^{\mu \nu }\varphi '_{x_\mu }\\
&& +\partial _{x_\nu
}(g^{ij})\varphi '_{x_i}\varphi '_{x_j}g^{\mu \nu }\varphi 
'_{x_\mu }+2g^{ij}\varphi
''_{x_i,x_\nu }\varphi '_{x_j}g^{\mu \nu }\varphi 
_{x_\mu }\nonumber\\
&=&
2\langle \varphi ''_{xx}\vert G\xi \otimes G\xi \rangle +2\langle \varphi
''_{xx}\vert G\varphi '_x\otimes G\varphi '_x\rangle +2\langle \partial _xG\vert
G\xi \otimes \varphi '_x\otimes \xi \rangle \nonumber\\
&& -\langle \partial
_xG\vert G\varphi '_x\otimes \xi \otimes \xi \rangle 
+\langle \partial
_xG\vert G\varphi '_x\otimes \varphi '_x\otimes \varphi 
'_x\rangle
.\nonumber
\end{eqnarray} 
Here we use the straight forward scalar products between tensors of the 
same size (2 ohr 3) and consider that the first index in the 3 tensor 
$\partial _xG$ is the one corresponding to the differentiations 
$\partial _{x_j}$. We also notice that $\varphi '_x,\xi $ are naturally 
cotangent vectors, while $G\varphi '_x,G\xi $ are tangent vectors. We want 
this 
Poisson bracket to be $\ge 0$ or even $\equiv 0$ on the set $a=b=0$, 
i.e. on the set given by 
\ekv{vc.6}
{
\langle G\vert \xi \otimes \xi -\varphi '_x\otimes \varphi '_x\rangle =0,\ 
\langle G\vert
\varphi '_x\otimes \xi \rangle =0.}

\par
\noindent \it Observation 1. \rm If $\varphi $
is a distance \fu{} in the sense that $\langle G\vert\varphi '_x\otimes \varphi 
'_x\rangle \equiv 1$, then if we differentiate in the direction 
$G\varphi '_x$, we get 
$$0=(G\varphi '_x\cdot \partial _x)\langle G\vert \varphi '_x\otimes \varphi 
'_x\rangle =\langle \partial _xG\vert G\varphi '_x\otimes \varphi '_x\otimes 
\varphi '_x\rangle +2\langle \varphi ''_{xx}\vert G\varphi '_x\otimes G\varphi 
'_x\rangle .$$
From this we see that two terms in the final expression in (\ref{vc.5}) 
cancel and we get 
\ekv{vc.7}
{{1\over 2}H_{a}b=2\langle \varphi ''_{xx}\vert G\xi \otimes G\xi \rangle 
+2\langle 
\partial _xG\vert G\xi \otimes \varphi '_x\otimes \xi \rangle -\langle 
\partial _xG\vert G\varphi '_x\otimes \xi \otimes \xi \rangle. }
\smallskip

\par
\noindent \it Observation 2. \rm If we replace 
$\varphi (x)$ by $\psi (x)=f(\varphi (x))$, then 
\begin{eqnarray*}
\psi '_x&=& f'(\varphi (x))\varphi '_x\\
\psi ''_{xx}&=& f''(\varphi (x))\varphi '_x\otimes \varphi '_x+f'(\varphi (x))\varphi 
''_{xx}. 
\end{eqnarray*}
If $\xi$  satisfies (\ref{vc.6}), then it is natural to replace $\xi $ by 
$\eta =f'(\varphi )\xi $, in order to preserve this condition (for the new 
symbol) and we see that all terms in the final member of (\ref{vc.5}), 
when restricted to $a=b=0$, become multiplied by $f'(\varphi )^3$ except 
the second one which becomes replaced by 
$$f'(\varphi )^32\langle \varphi ''_{xx}\vert G\varphi '_x\otimes G\varphi '_x
\rangle +2f''(\varphi (x))f'(\varphi (x))^2\langle G\vert \varphi '_x\otimes  
\varphi '_x\rangle ^2.$$
(For the first term in (\ref{vc.5}) we also use that $\langle \varphi '_x 
\otimes \varphi '_x\vert G\xi \otimes G\xi \rangle =\langle \varphi '_x\vert G 
\xi \rangle ^2=0.$) Thus we get after the two substitutions $\varphi \mapsto 
\psi =f(\varphi (x))$, $\xi \mapsto \eta =f'(\varphi (x))\xi $:
\begin{eqnarray}\label{vc.8}
{1\over 2}H_{a}b(x,\eta )&=& 2f''(\varphi (x))f'(\varphi (x))^2\Vert \varphi 
'_x\Vert 
_g^4+\\
&& f'(\varphi )^3\big( 2\langle \varphi ''_{xx}\vert G\xi \otimes G\xi \rangle 
+2\langle \varphi ''_{xx}\vert G\varphi '_x\otimes G\varphi '_x\rangle +2\langle 
\partial _xG\vert G\xi \otimes \varphi '_x\otimes \xi \rangle \nonumber\\
&& -\langle \partial _xG\vert G\varphi '_x\otimes \xi \otimes \xi \rangle 
+\langle \partial _xG\vert G\varphi '_x\otimes \varphi '_x\otimes \varphi 
'_x\rangle \big) , \nonumber
\end{eqnarray}
with $\eta =f'(\varphi )\xi $, $\xi $ satisfying (\ref{vc.6}), so that 
$\eta $ satisfies the same condition (with $\varphi $ replaced by $\psi $):
\ekv{vc.9}
{
\langle G\vert \eta \otimes \eta -\psi '_x\otimes \psi '_x\rangle =
\langle G\vert\psi '_x\otimes \eta \rangle =0.}
Moreover $\Vert \varphi '_x\Vert ^2_g=\langle G\vert \varphi '_x\otimes \varphi 
'_x\rangle $ by definition.
\smallskip

\par
\noindent \it Conclusion. \rm To get $H_{a}b\ge 0$ whenever 
(\ref{vc.6}) is satisfied, it suffices to start with a \fu{} $\varphi $ 
with non-vanishing gradient, and then replace $\varphi $ by $f(\varphi )$ with 
$f'>0$ and $f''/f'$ \sufly{} large. This kind of convexification ideas are 
very old and used recently in a related context by Lebeau--Robbiano 
\cite{LeRo}, Burq \cite{Bu}. For later use, we needed to spell out the 
calculations quite explicitly.

\section{Carleman estimate}\label{SectionCe}
\setcounter{equation}{0}

We use from now  on semiclassical notation (see for instance \cite{DiSj}).

Let $P_0=-h^2\Delta =\sum (hD_{x_j})^2$, with $D_{x_j}={1\over i}\partial 
_{x_j}$. Let $\varphi$, $\widetilde{\Omega }$ be as in the beginning of 
Section \ref{SectionVC}. Then 
\ekv{ce.1}
{e^{\varphi /h}\circ P_0\circ e^{-\varphi /h}=\sum_{j=1}^n(hD_{x_j}+i\partial 
_{x_j}\varphi )^2=A+iB,}
where $A,B$ are the formally \sa{} \op{}s:
\ekv{ce.2}
{A=(hD)^2-(\varphi '_x)^2,\ B=\sum(\partial _{x_j}\varphi \circ hD_{x_j} 
+ hD_{x_j}\circ \partial _{x_j}\varphi )}
having the Weyl symbols (for the semi-classical quantization)
\ekv{ce.3}
{a=\xi ^2-(\varphi '_x)^2,\ b=2\varphi '_x\cdot \xi .}

\par We assume that $\varphi $ has non vanishing gradient
and is a limiting Carleman weight 
in the sense that 
\ekv{ce.4}
{
\{a,b\} (x,\xi )=0,\hbox{ when }a(x,\xi )=b(x,\xi
)=0.
}
Here $\{ a,b\}=a_\xi '\cdot b'_x-a'_x\cdot b'_\xi $ is the Poisson
bracket (as in (\ref{vc.5})):
\ekv{ce.5}
{\{ a,b\} =4\langle \varphi ''_{xx}(x)\vert\xi \otimes \xi +\varphi '_x\otimes 
\varphi 
'_x\rangle .
}
On the $x$-dependent hypersurface in $\xi $-space, given by $b(x,\xi )=0$, 
we know that the quadratic \pol{} $\{ a,b\}(x,\xi )$ vanishes when $\xi 
^2=(\varphi '_x)^2$. It follows that 
\ekv{ce.6}
{
\{ a,b\} (x,\xi )=c(x)(\xi ^2-(\varphi '_x)^2),\hbox{ for }b(x,\xi )=0,
}
where $c(x)\in C^\infty (\widetilde{\Omega };{\bf R})$. Then consider
$$\{ a,b\} (x,\xi )-c(x)(\xi ^2-(\varphi '_x)^2),$$
which is a quadratic \pol{} in $\xi $, vanishing when $\varphi '_x(x)\cdot 
\xi =0$ It follows that this is of the form $\ell (x,\xi )b(x,\xi )$ 
where $\ell (x,\xi )$ is affine in $\xi $ with smooth coefficients, and we 
end up with
\ekv{ce.7}
{
\{ a,b\} =c(x)a(x,\xi )+\ell (x,\xi )b(x,\xi ).
}
But $\{ a,b\}$ contains no linear terms in $\xi $, so we know that $\ell 
(x,\xi )$ is linear in $\xi $.
\par The commutator $[A,B]$ can be computed directly: and we get
%$$
%[A,B]=-{h\over i}[(h\partial )^2+(\varphi '_x)^2,\varphi '(x)\circ \partial 
%_x+\partial _x\circ \varphi '(x)].
%$$
%Here 
%\begin{eqnarray*}
%&[(h\partial )^2,\varphi '(x)\circ \partial _x]=\sum_{j,k}[(h\partial 
%_{x_j})^2,\varphi '_{x_k}(x)\circ \partial _{x_k}]=&\\
%&=\sum_{j,k}[(h\partial _{x_j})^2,\varphi _{x_k}'(x)]\partial 
%_{x_k}=\sum_{j,k} (h\partial _{x_j}\circ \varphi ''_{x_jx_k}h\partial 
%{x_k}+\varphi ''_{x_jx_k}h\partial _{x_j}h\partial 
%_{x_k})=&\\
%&=\sum_{j,k}(h\partial _{x_j}\circ \varphi ''_{x_jx_k}+\varphi 
%''_{x_jx_k}\circ h\partial _{x_j})\circ h\partial _{x_k},&
%\end{eqnarray*}
%$$
%[(h\partial )^2,\partial _x\circ \varphi '(x)]=\sum_{j,k}h\partial 
%_{x_k}\circ (h\partial _{x_j}\circ \varphi ''_{x_jx_k}+\varphi ''_{x_jx_k}\circ 
%h\partial _{x_j}),
%$$
%$$
%[(\varphi '_x)^2,\varphi '(x)\circ \partial _x]=-\varphi '(x)\cdot \partial 
%_x((\varphi '_x)^2)=-2\varphi ''_{xx}(x)\varphi '(x)\cdot \varphi '(x).
%$$
%$$
%[(\varphi '_x)^2,\partial _x\circ \varphi '(x)]=-2\varphi 
%''_{xx}(x)\varphi '(x)\cdot \varphi '_x(x).
%$$
%Thus,
\begin{eqnarray*}
[A,B]&=&
%-{h\over i}\big( \sum_{j,k}[(h\partial _{x_j}\circ \varphi 
%''_{x_jx_k}+\varphi ''_{x_jx_k}\circ h\partial _{x_j})h
%\partial _{x_k}+h\partial 
%_{x_k}\circ (h\partial _{x_j}\circ \varphi ''_{x_jx_k}+\varphi ''_{x_jx_k}
%\circ h\partial _{x_j})]
%\\
%&&\hskip 3cm -4\langle \varphi ''_{xx};\varphi '(x)\otimes \varphi '(x)\rangle 
%\big)\\
{h\over i}\Big( \sum_{j,k}[(hD_{x_j}\circ \varphi ''_{x_jx_k}+\varphi 
''_{x_jx_k}\circ hD_{x_j})hD_{x_k}+
hD_{x_k}(hD_{x_j}\circ \varphi ''_{x_jx_k}+
\varphi ''_{x_jx_k}\circ hD_{x_j})]\\
&&\hskip 3cm +4\langle \varphi ''_{xx},\varphi '_x(x)\otimes 
\varphi '_x(x)\rangle \Big).
\end{eqnarray*}

\par The Weyl symbol of $[A,B]$ as a semi-classical \op{} is 
$${h\over i}\{ a,b\} +h^3p_0(x),$$
 Combining this with (\ref{ce.7}), we get with a new  
$p_0$:
\ekv{ce.8}
{
i[A,B]=h({1\over 2}(c(x)\circ A+A\circ c)+
{1\over 2}(LB+BL)+h^2p_0(x)),
}
where $L$ denotes the Weyl quantization of $\ell$.

\def\Ca{Carleman}

\par We next derive the \Ca{} estimate for $u\in C_0^\infty (\Omega )$, 
$\Omega 
\subset\subset \widetilde{\Omega }$:
Start from $P_0u=v$ and let $\widetilde{u}=e^{\varphi /h}u$, 
$\widetilde{v}=e^{\varphi /h}v$, so that
\ekv{ce.9}
{
(A+iB)\widetilde{u}=\widetilde{v}.
}
Using the formal \sa{}ness of $A,B$, we get 
\ekv{ce.10}
{
\Vert \widetilde{v}\Vert ^2=((A-iB)(A+iB)\widetilde{u}\vert 
\widetilde{u})=\Vert A\widetilde{u}\Vert ^2+\Vert B\widetilde{u}\Vert 
^2+(i[A,B]\widetilde{u}\vert \widetilde{u}).
}
Using (\ref{ce.8}), we get for $u\in C_0^\infty (\Omega )$:
\begin{eqnarray}\label{ce.11}
\Vert \widetilde{v}\Vert ^2&\ge&\Vert A\widetilde{u}\Vert ^2+\Vert 
B\widetilde{u}\Vert ^2-{\cal O}(h)(\Vert A\widetilde{u}\Vert \Vert 
\widetilde{u}\Vert 
+\Vert L\widetilde{u}\Vert \Vert B\widetilde{u}\Vert )-
{\cal O}(h^3)\Vert \widetilde{u}\Vert ^2\\
&\ge & {2\over 3}\Vert A\widetilde{u}\Vert ^2+{1\over 2}\Vert 
B\widetilde{u}\Vert ^2-{\cal O}(h^2)(\Vert \widetilde{u}\Vert ^2+\Vert 
L\widetilde{u}\Vert ^2)
.\nonumber \\
&\ge &{1\over 2}(\Vert A\widetilde{u}\Vert ^2+\Vert B\widetilde{u}\Vert 
^2)-{\cal O}(h^2)\Vert \widetilde{u}\Vert ^2,\nonumber
\end{eqnarray}
where in the last step we used the apriori estimate 
$$
\Vert h\nabla \widetilde{u}\Vert ^2\le {\cal O}(1)(\Vert 
A\widetilde{u}\Vert ^2+\Vert \widetilde{u}\Vert ^2),
$$
which follows from the classical ellipticity of $A$.

Now we could try to use that $B$ is associated to a non-vanishing 
gradient field (and 
hence without any closed or even trapped trajectories in 
$\widetilde{\Omega }$), to obtain the Poincar{\'e} estimate:
\ekv{ce.12}
{h\Vert \widetilde{u}\Vert \le {\cal O}(1)\Vert B\widetilde{u}\Vert .}

\par We see that (\ref{ce.12}) is not quite good enough to absorb the last 
term in (\ref{ce.11}). In order to remedy for this, we make a slight 
modification of $\varphi $ by introducing
\ekv{ce.13}
{
\varphi _\epsilon =f\circ \varphi ,\hbox{ with }f=f_\epsilon 
}
to be chosen below, and write $a_\epsilon +ib_\epsilon $ for the 
conjugated 
symbol. We saw in Section \ref{SectionVC} and especially in (\ref{vc.8}) 
that the Poisson 
bracket $\{ a_\epsilon ,b_\epsilon \}$, becomes with $\varphi $ equal to the 
original weight:
\ekv{ce.14}
{\{ a_\epsilon ,b_\epsilon \} (x,f'(\varphi )\eta )=f'(\varphi )^3(\{ a,b\} 
(x,\eta )+{4f''(\varphi )\over f'(\varphi )}\Vert \varphi '_x\Vert ^4 ),\hbox{ when 
}a(x,\eta)=b(x,\eta)=0.}
The substitution $\xi \to f'(\varphi )\eta $ is motivated be the fact that if 
$a(x,\eta )=b(x,\eta )=0$, then $a_\epsilon (x,f'(\varphi )\eta )=b_\epsilon 
(x,f'(\varphi )\eta )=0$. Now let 
\ekv{ce.15}
{
f_\epsilon (\lambda )=\lambda +\epsilon \lambda ^2/2,
}
with 
$0\le \epsilon \ll 1$, so that 
$${4f''(\varphi )\over f'(\varphi )}={4\epsilon \over 1+\epsilon \varphi }= 
4\epsilon +{\cal O}(\epsilon ^2).$$
In view of (\ref{ce.14}), (\ref{ce.4}), we get 
\ekv{ce.16}
{
\{ a_\epsilon ,b_\epsilon \} (x,\xi )=4f''_\epsilon (\varphi )(f'_\epsilon 
(\varphi ))^2\Vert \varphi '\Vert ^4\approx 4\epsilon \Vert \varphi '_x\Vert ^4,
}
when $a_\epsilon (x,\xi )=b_\epsilon (x,\xi )=0$, so instead of 
(\ref{ce.7}), we get 
\ekv{ce.17}
{
\{a_\epsilon ,b_\epsilon \} =4f''_\epsilon (\varphi )(f'_\epsilon (\varphi 
))^2\Vert \varphi '_x\Vert ^4+c_\epsilon (x)a_\epsilon (x,\xi )+\ell 
_\epsilon 
(x,\xi )b_\epsilon (x,\xi ),
}
with $\ell_\epsilon (x,\xi )$ linear in $\xi $.
%\ekv{ce.18}
%{
%i[A_\epsilon ,B_\epsilon ]=h2f''_\epsilon (\varphi )(f'_\epsilon (\varphi 
%))^2\Vert \varphi '_x\Vert ^4+h({1\over 2}(c_\epsilon (x)\circ A_\epsilon 
%+A_\epsilon \circ c_\epsilon )+B_\epsilon \circ d_\epsilon \circ 
%B_\epsilon +h^2p_2(x)).}

\par Instead of (\ref{ce.11}), we get with $\widehat{u}=e^{\varphi _\epsilon 
/h}u$, $\widehat{v}=e^{\varphi _\epsilon /h}v$ when $P_0u=v$:
\ekv{ce.19}
{
\Vert \widehat{v}\Vert ^2\ge h(4\epsilon +{\cal O}(\epsilon ^2))\int 
\Vert \varphi '_x\Vert ^4\vert \widehat{u}(x)\vert ^2dx +{1\over 2}\Vert 
A_\epsilon \widehat{u}\Vert ^2+{1\over 2}\Vert B_\epsilon 
\widehat{u}\Vert ^2-{\cal O}(h^2)\Vert \widehat{u}\Vert ^2, 
} while the analogue of (\ref{ce.12}) remains \ufly{} valid when 
$\epsilon $ is small:
\ekv{ce.20}
{h\Vert \widehat{u}\Vert \le {\cal O}(1)\Vert B_\epsilon \widehat{u}\Vert 
,}
even though we will not use this estimate.

\par Choose $h\ll\epsilon \ll 1$, so that (\ref{ce.19}) gives
\ekv{ce.21}
{
\Vert \widehat{v}\Vert ^2\ge \epsilon h\Vert \widehat{u}\Vert 
^2+{1\over 2}\Vert A_\epsilon \widehat{u}\Vert ^2+{1\over 2}\Vert 
B_\epsilon \widehat{u}\Vert ^2.
} 

\par We want to transform this into an estimate for 
$\widetilde{u},\widetilde{v}$. From the special form of $A_\epsilon $, we 
see that 
$$\Vert hD\widehat{u}\Vert ^2\le (A_\epsilon \widehat{u}\vert 
\widehat{u})+{\cal O}(1)\Vert \widehat{u}\Vert ^2,$$
leading to 
$$
\Vert hD\widehat{u}\Vert ^2 \le {1\over 2}\Vert A_\epsilon 
\widehat{u}\Vert 
^2+{\cal O}(1)\Vert \widehat{u}\Vert ^2.
$$
Combining this with (\ref{ce.21}), we get 
\ekv{ce.22}
{
\Vert \widehat{v}\Vert ^2\ge {\epsilon h\over C_0}(\Vert \widehat{u}\Vert 
^2+\Vert hD\widehat{u}\Vert ^2)+({1\over 2}-{\cal O}(\epsilon h))\Vert 
A_\epsilon \widehat{u}\Vert ^2+{1\over 2}\Vert B_\epsilon 
\widehat{u}\Vert ^2.
}
Write $\varphi _\epsilon =\varphi +\epsilon g$, where $g=g_\epsilon $ is ${\cal 
O}(1)$ with all its derivatives. We have 
$$\widehat{u}=e^{\epsilon g/h}\widetilde{u},\ 
\widehat{v}=e^{\epsilon g/h}\widetilde{v},$$
so 
$$hD\widehat{u}=e^{\epsilon g /h}(hD\widetilde{u}+{\epsilon \over i} 
g'\widetilde{u})=e^{{\epsilon \over h}g}(hD\widetilde{u}+{\cal 
O}(\epsilon )\widetilde{u}), $$
and 
\begin{eqnarray*}
\Vert \widehat{u}\Vert ^2+\Vert hD\widehat{u}\Vert ^2 &\ge & \Vert 
e^{\epsilon g/h}\widetilde{u}\Vert ^2+\Vert e^{\epsilon 
g/h}hD\widetilde{u}\Vert ^2 -C\epsilon \Vert e^{\epsilon g/h }
\widetilde{u}\Vert \Vert e^{\epsilon g/h}hD\widetilde{u}\Vert -C\epsilon 
^2\Vert e^{\epsilon g/h}\widetilde{u}\Vert ^2\\
&\ge & (1-C\epsilon )(\Vert e^{\epsilon g/h}\widetilde{u}\Vert ^2+\Vert 
e^{\epsilon g/h}hD\widetilde{u}\Vert ^2),
\end{eqnarray*}
so from (\ref{ce.22}) we obtain after increasing $C_0$ by a factor 
$(1+{\cal O}(\epsilon ))$:
\ekv{ce.23}
{\Vert e^{\epsilon g/h}\widetilde{v}\Vert ^2\ge {\epsilon h\over 
C_0}(\Vert 
e^{\epsilon g/h}\widetilde{u}\Vert ^2+\Vert e^{\epsilon 
g/h}hD\widetilde{u}\Vert ^2).}

\par If we take $\epsilon =Ch$ with $C\gg 1$ but fixed, then $\epsilon 
g/h$ is \ufly{} \bdd{} in $\Omega $ and we get the Carleman estimate
\ekv{ce.24}
{
h^2(\Vert \widetilde{u}\Vert ^2+\Vert hD\widetilde{u}\Vert ^2)\le C_1 
\Vert \widetilde{v}\Vert ^2.}
This clearly extends to solutions of the equation 
\ekv{ce.24.5}{(-h^2\Delta +h^2q)u=v,}
if $q\in L^\infty $ is fixed, since we can start by applying
(\ref{ce.24}) with $\widetilde{v}$ replaced by 
$\widetilde{v}-h^2q\widetilde{u}$. 
Summing up the discussion so far, we have
\begin{prop}\label{PropCE.1}
Let $P_0$, $\widetilde{\Omega }$, $\varphi $ be as in the beginning of this 
section and assume that $\varphi $ is a limiting Carleman weight in the sense 
that (\ref{ce.4}) holds. Let $\Omega \subset\subset \widetilde{\Omega }$ 
be open and let $q\in L^\infty (\Omega )$. Then if $u\in C_0^\infty 
(\Omega )$, we have 
\ekv{ce.24.6}
{
h(\Vert e^{\varphi /h}u\Vert +\Vert hDe^{\varphi /h}u\Vert )\le C \Vert e^{\varphi 
/h}(-h^2\Delta +h^2q)u\Vert ,
}
where $C$ depends on $\Omega $, and $h>0$ is small enough so that $Ch\Vert 
q\Vert _{L^\infty (\Omega )}\le 1/2$.\end{prop}

\par We next establish a Carleman estimate when $P_0u=v$, $u\in C^\infty 
(\Omega )$, ${u_\vert}_{\partial \Omega }=0$ and $\Omega \subset\subset 
\widetilde{\Omega }$ is a domain with $C^\infty $ \bdy{}. As before, we 
let $\widehat{u}=e^{\varphi /h}u$, $\widehat{v}=e^{\varphi /h}v$, with 
$\varphi =\varphi _\epsilon $, $0\le \epsilon \ll 1$. With $A=A_\epsilon $, 
$B=B_\epsilon $, we have 
\ekv{ce.25}
{(A+iB)\widehat{u}=\widehat{v},}
and 
\begin{eqnarray}\label{ce.26}
\Vert \widehat{v}\Vert ^2 &=& ((A+iB)\widehat{u}\vert (A+iB)\widehat{u})\\
&=& \Vert A\widehat{u}\Vert ^2+\Vert B\widehat{u}\Vert 
^2+i((B\widehat{u}\vert A\widehat{u})-(A\widehat{u}\vert B\widehat{u})),
\nonumber
\end{eqnarray}
Using that $B$ is a first order \dop{} and that 
$${\widehat{u}\big\vert }_{\partial \Omega }=0,$$ we see that 
\ekv{ce.27}
{
(A\widehat{u}\vert B\widehat{u})=(BA\widehat{u}\vert \widehat{u}).
}
Similarly, we have 
\ekv{ce.28}
{
(B\widehat{u}\vert (\varphi '_x)^2\widehat{u})=((\varphi 
'_x)^2B\widehat{u}\vert \widehat{u}).
}
Finally, we use Green's formula, with $\nu $ denoting the exterior unit 
normal, to transform 
$$
(B\widehat{u}\vert -h^2\Delta \widehat{u})_\Omega =-h^2(B\widehat{u}\vert 
\partial _\nu \widehat{u})_{\partial \Omega }+(-h^2\Delta 
B\widehat{u}\vert \widehat{u})_\Omega ,
$$
where we also used that ${\widehat{u}\big\vert}_{\partial \Omega }=0$.

\par On $\partial \Omega $, we have 
$$B=2(\varphi '_x\cdot \nu ){h\over i}\partial _\nu +B',$$
where $B'$ acts along the \bdy{}, so using again the Dirichlet condition, 
we get 
$$(B\widehat{u}\vert \partial _\nu \widehat{u})_{\partial \Omega 
}={2h\over 
i}((\varphi '_x\cdot \nu )\partial _\nu \widehat{u}\vert \partial _\nu 
\widehat{u})_{\partial \Omega }.$$

Putting together the calculations and using (\ref{ce.2}) for $A$, we get
\ekv{ce.29}
{\Vert \widehat{v}\Vert^2 =\Vert A\widehat{u}\Vert ^2+\Vert 
B\widehat{u}\Vert ^2+i([A,B]\widehat{u}\vert \widehat{u})-2h^3((\varphi 
'_x\cdot \nu )\partial _\nu \widehat{u}\vert \partial _\nu 
\widehat{u})_{\partial \Omega }.}

\par Let
\ekv{ce.29.5}
{
\partial \Omega _{\pm}=\{ x\in\partial \Omega ;\, \pm \varphi '_x\cdot  \nu 
\ge 0\}.
}
Notice that $\partial \Omega _\pm$
are \indep{} of $\epsilon $. We rewrite (\ref{ce.29}) as 
\ekv{ce.30}
{
-2h^3((\varphi '_x\cdot \nu )\partial _\nu \widehat{u}\vert \partial _\nu 
\widehat{u})_{\partial \Omega _-}+i([A,B]\widehat{u}\vert 
\widehat{u})+\Vert 
A\widehat{u}\Vert ^2+\Vert B\widehat{u}\Vert ^2=\Vert \widehat{v}\Vert 
^2+2h^3((\varphi '_x\cdot \nu )\partial _\nu \widehat{u}\vert \partial _\nu 
\widehat{u})_{\partial \Omega _+}.}

\par This is analogous to (\ref{ce.10}) and the extra \bdy{} terms can be 
added in the discussion leading from (\ref{ce.19}) to (\ref{ce.22}) and we 
get instead of (\ref{ce.22}):
\begin{eqnarray}\label{ce.31}
&-2h^3((\varphi '_x\cdot \nu )\partial _\nu \widehat{u}\vert \partial _\nu 
\widehat{u})_{\partial \Omega _-}+{\epsilon h\over C_0}(\Vert 
\widehat{u}\Vert ^2+\Vert hD\widehat{u}\Vert ^2)+&\\
&({1\over 2}-{\cal 
O}(\epsilon h))\Vert A_\epsilon \widehat{u}\Vert ^2+{1\over 2}\Vert 
B_\epsilon \widehat{u}\Vert ^2&\nonumber\\
&\hskip 2cm \le  
\Vert \widehat{v}\Vert ^2+2h^3((\varphi '_x\cdot \nu )\partial _\nu 
\widehat{u}\vert \partial _\nu \widehat{u})_{\partial \Omega _+},&\nonumber
\end{eqnarray}
with $\varphi =\varphi _\epsilon $, provided $\epsilon \gg h$. Fixing $\epsilon 
=Ch$ for $C\gg 1$, we get with $\varphi =\varphi _{\epsilon =0}$ for some 
$C_0>0$:
\begin{eqnarray}\label{ce.32}
&-{h^3\over C_0}((\varphi '_x\cdot \nu )\partial _\nu \widetilde{u}\vert 
\partial _\nu \widetilde{u})_{\partial \Omega _-}+{h^2\over C_0}(\Vert 
\widetilde{u}\Vert ^2+\Vert hD\widetilde{u}\Vert ^2)&\\
&
\le \Vert \widetilde{v}\Vert ^2+C_0h^3((\varphi '_x\cdot \nu 
)\partial _\nu \widetilde{u}\vert \partial _\nu \widetilde{u})_{\partial 
\Omega _+}.
&\nonumber
\end{eqnarray}
Here we recall that $-h^2\Delta u=v$, $\widetilde{u}=e^{\varphi /h}u$, 
$\widetilde{v}=e^{\varphi /h}v$, $\varphi =\varphi _{\epsilon =0}$, 
${u_\vert}_{\partial \Omega }=0$. 

\par If $q\in L^\infty $, we get for $h^2(-\Delta +q)u=v$, 
${u_\vert}_{\partial \Omega }=0$, by applying (\ref{ce.32}) with 
$\widetilde{v}$ replaced by $\widetilde{v}-h^2q\widetilde{u}$:
\begin{eqnarray}
\label{ce.33}
&-{h^3\over C_0}((\varphi '_x\cdot \nu )\partial _\nu \widetilde{u}\vert 
\partial _\nu \widetilde{u})_{\partial \Omega _-}+{h^2\over C_0}(\Vert 
\widetilde{u}\Vert ^2+\Vert hD\widetilde{u}\Vert ^2)&\\
&
\le \Vert \widetilde{v}\Vert ^2+C_0h^3((\varphi '_x\cdot \nu 
)\partial _\nu \widetilde{u}\vert \partial _\nu \widetilde{u})_{\partial 
\Omega _+}.
&\nonumber
\end{eqnarray}
Here $\widetilde{u},\widetilde{v}$
are defined as before.

\par Summing up, we have
\begin{prop}\label{PropCE.2}
Let $\widetilde{\Omega },\varphi $ be as in Proposition \ref{PropCE.1}. 
Let $\Omega \subset\subset \widetilde{\Omega }$ be an open set with 
$C^\infty $ \bdy{} and let $q\in L^\infty (\Omega )$. Let $\nu $ denote 
the exterior unit normal to $\partial \Omega $ and define $\partial \Omega 
_{\pm}$ as in (\ref{ce.29.5}). Then there exists a constant $C_0>0$, such 
that for every $u\in C^\infty (\overline{\Omega }) $ with 
${u_\vert}_{\partial \Omega }=0$, we have for $0<h\ll 1$:
\begin{eqnarray}\label{ce.24.7}
&
-{h^3\over C_0}((\varphi '_x\cdot \nu )e^{\varphi /h}\partial _\nu u\vert 
e^{\varphi /h}\partial _\nu u)_{\partial \Omega _-}+{h^2\over C_0}(\Vert 
e^{\varphi /h}u\Vert ^2+\Vert e^{\varphi /h}h\nabla u\Vert ^2)
&\\
&
\le \Vert e^{\varphi /h}(-h^2\Delta +h^2q)u\Vert ^2+C_0h^3((\varphi '_x\cdot 
\nu )e^{\varphi /h}\partial _\nu u\vert e^{\varphi /h}\partial _\nu 
u)_{\partial \Omega _+},
&\nonumber
\end{eqnarray}
\end{prop}
\medskip

\par\noindent \it Remark. \rm If $\varphi $ is a limiting Carleman weight, 
then so is $-\varphi $. With $\widetilde{u}=e^{-\varphi /h}u$, 
$\widetilde{v}=e^{-\varphi /h}v$, we still have (\ref{ce.33}), provided we 
permute $\partial \Omega _-$
and $\partial \Omega _+$
and change the signs in front of the \bdy{} terms, so that they remain 
positive. 
\medskip

\section{Construction of complex geometrical optics 
solutions}\label{SectionSP}
\setcounter{equation}{0}

\par Let $H^s({\bf R}^n)$ denote the semi-classical Sobolev space of 
order $s$, equipped with the norm $\Vert \langle hD\rangle ^su\Vert $. We 
define $H^s(\Omega )$, $H_0^s(\Omega )$ in the usual way, when $\Omega 
\subset\subset {\bf R}^n$ has smooth \bdy{}. (\ref{ce.24}) can be written 
\ekv{sp.1}
{
h\Vert {u}\Vert_{H^1} \le C\Vert e^{\varphi /h}P_0e^{-\varphi 
/h}{u}\Vert ,\ {u}\in C_0^\infty (\Omega ),
}
when $P_0=-h^2\Delta$. Here we let $\Omega \subset \widetilde{\Omega }$ be 
as in Section \ref{SectionCe}. Recall that $P_{0,\varphi }=e^{\varphi /h}P_0 
e^{-\varphi /h}$ has the semiclassical Weyl symbol $\xi ^2-{\varphi 
'_x}^2+2i\varphi 
'_x\cdot \xi =a+ib$, which is elliptic in the region $\vert \xi \vert \ge 
2\vert \varphi '(x)\vert $. It is therefore clear that (\ref{sp.1}) can be 
extended to: 
\ekv{sp.2}
{
h\Vert u\Vert _{H^{-s+1}}\le C_{s,\Omega }\Vert e^{\varphi /h}P_0e^{-\varphi 
/h}u\Vert _{H^{-s}},\ u\in C_0^\infty (\Omega ),
}
for every fixed $s\in {\bf R}$. With $q\in L^\infty (\widetilde{\Omega 
})$, 
we put 
$$P=-h^2(\Delta -q),\ P_\varphi =e^{\varphi /h}Pe^{-\varphi /h}=P_{0,\varphi }+h^2q.$$

\par If $0\le s\le 1$, we have $$\Vert qu\Vert _{H^{-s}}\le \Vert qu\Vert 
\le 
\Vert q\Vert _{L^\infty }\Vert 
u\Vert \le \Vert q\Vert _{L^\infty }\Vert u\Vert _{H^{-s+1}},
$$
and for $h>0$ small enough, we get from (\ref{sp.2}):
\ekv{sp.3}
{
h\Vert u\Vert _{H^{-s+1}}\le C_{s,\Omega }\Vert e^{\varphi /h}Pe^{-\varphi /h} 
u\Vert _{H^{-s}}.
}
The Hahn-Banach theorem now implies in the usual way:
\begin{prop}\label{PropSP.1}
Let $0\le s\le 1$. Then for $h\ge 0$ small enough, for every 
$v\in H^{s-1}(\Omega )$, there exists $u\in H^s(\Omega )$ such that 
\ekv{sp.4}
{
e^{-\varphi /h}Pe^{\varphi /h}u=v,\ h\Vert u\Vert _{H^s}\le C\Vert v\Vert 
_{H^{s-1}}.
}
\end{prop}

\par This result remains valid, when $q$ is complex valued. In that case 
we 
replace $P$ in (\ref{sp.3}) by $\overline{P}=-h^2\Delta +\overline{q}$. 

\par We next construct certain WKB-solutions to the homogeneous equation. 
Recall that $a,b$ are in involution on the joint zero set $J:\, a=b=0$ in 
view of (\ref{ce.7}). At the points of $J$ we also see that the \hf{}s 
\ekv{sp.5}
{
H_a=2(\xi \cdot \partial _x+\langle \varphi ''_{xx}\varphi '_x\vert \partial 
_\xi\rangle  ),\ 
H_b=2(\varphi '_x\cdot \partial _x-\langle \varphi ''_{xx}\xi \vert \partial _\xi
\rangle  )
}
are linearly \indep{} and even have linearly \indep{} $x$-space 
projections. We conclude that $J$ is an involutive \mfld{} such that each 
bicharacteristic leaf (of dimension 2) has a base space projection which 
is also a nice sub\mfld{} of dimension 2. It follows that we have plenty 
of 
smooth local solutions to the Hamilton-Jacobi problem
\ekv{sp.6}
{
a(x,\psi '(x))=b(x,\psi '(x))=0.
} 
Indeed, if $(x_0,\xi _0)\in J$, and we let $H\subset \Omega $ be a 
sub\mfld{} of codimension $2$ passing through $x_0$ transversally to the 
projection 
of the bicharacteristic leaf through $(x_0,\xi _0)$, then we have a unique 
local solution of (\ref{sp.6}), with ${\psi 
_\vert}_{H}=\widetilde{\psi }$, 
if $\widetilde{\psi }$ is a smooth 
real-valued \fu{} on $H$ such that $\widetilde{\psi }'(x_0)$ is equal to 
the projection of $(x_0,\xi _0)$ in $T^*_{x_0}(H)$.

\par Since we need some explicit control of the size of the domain of 
definition of $\psi $, we now give a more down-to-earth construction. 
(\ref{sp.6}) can be written more explicitly as 
\ekv{sp.7}
{
\psi '(x)^2-\varphi '(x)^2=0,\ \varphi '(x)\cdot \psi '(x)=0.
}
First restrict the attention to the hypersurface $G=\varphi ^{-1}(C_0)$ for 
some fixed constant $C_0$, and let $g$ denote the restriction of $\psi $ 
to $G$. Then we get the necessary condition that 
\ekv{sp.8}
{
g'(x)^2=\varphi '(x)^2,}
where $g'(x)^2$ is the square of the norm of the differential for the 
metric 
dual to $e_0$, the induced Euclidean metric. Now (\ref{sp.8}) is a 
standard eikonal \e{} on $G$ and we can find solutions of the form 
$g(x)={\rm dist\,}(x,\Gamma )$, where $\Gamma $ is either a point or a 
hypersurface in $G$ and ${\rm dist}$ denotes the distance on $G$ with 
respect to the metric $\varphi '(x)^2e_0(dx)$. Of course, we will have to be 
careful, since such distance functions in general will develop 
singularities, and in the following we restrict $G$ if necessary, so that 
the function $g$ is smooth. With $g$ solving (\ref{sp.8}), we define $\psi 
$ 
to be the extension of $g$ which is constant along the integral curves of 
the field $\varphi '(x)\cdot \partial _x$:
\ekv{sp.9}
{
\varphi '(x)\cdot \partial _x\psi (x)=0,\ {\psi _\vert}_{G}=g.
}
Then the second equation in (\ref{sp.7}) holds by construction, and the 
first \e{} is fulfilled at the points of $G$. In order to verify that \e{} 
also away from $G$, we consider,
\ekv{sp.10}
{
\varphi '(x)\cdot \partial _x({\psi '}^2-{\varphi '}^2)=2(\langle \psi ''\varphi 
'\vert \psi '\rangle -\langle \varphi ''\varphi '\vert \varphi '\rangle ).
}
Taking the gradient of the second \e{} in (\ref{sp.7}), we get $\varphi 
''\psi '+\psi ''\varphi '=0$, and hence 
\begin{eqnarray}
\varphi '(x)\cdot \partial _x({\psi '}^2-{\varphi '}^2)&=&-2(\langle \varphi 
''\psi '\vert \psi '\rangle +\langle \varphi ''\varphi '\vert \varphi '\rangle 
)=-{1\over 2}\{ a,b\} (x,\psi ')\\
&=& -{1\over 2}c(x)({\psi '}^2-{\varphi '}^2)-\ell(x,\psi ')\varphi '\cdot \psi 
'\nonumber\\
&=& -{1\over 2}c(x)({\psi '}^2-{\varphi '}^2).
\end{eqnarray} 
Thus 
\ekv{sp.12}
{
(\varphi '(x)\cdot \partial _x+{c(x)\over 2})({\psi '}^2-{\varphi '}^2)=0,\ 
{({\psi '}^2-{\varphi '}^2)_\vert }_G=0,}
and we conclude that ${\psi '}^2-{\varphi '(x)}^2=0$.
\smallskip

\par \it Summing up the discussion so far, we have seen that if $\varphi $ is 
a 
limiting Carleman weight, and the open set $\Omega $
is a union of integral segments of $\varphi '(x)\cdot \partial _x$ all 
crossing the smooth hypersurface $G\subset \varphi ^{-1}(C_0) $, then if $g$
is smooth solution to the eikonal equation (\ref{sp.8}) on $G$ and we 
define $\psi $ to be the solution of (\ref{sp.9}), we get a solution 
of (\ref{sp.6}).\rm
\smallskip

\par (\ref{sp.6}) implies that 
\ekv{sp.13}
{
p(x,i\varphi '(x)+\psi '(x))=0,
}
which is the eikonal equation for the construction of WKB-solutions of the 
form $u(x;h)=a(x;h)e^{{1\over h}(-\varphi +i\psi )}$ of $P_0u\approx 0$. If 
we try $a$ smooth and \indep{} of $h$, we get
\begin{eqnarray}\label{sp.14}
e^{-{1\over h}(-\varphi +i\psi )}P_0e^{{1\over h}(-\varphi +i\psi 
)}a&=&e^{-{i\over 
h}\psi }P_{0,\varphi }e^{{i\over h}\psi }a\\
&=&\big( ((hD+\psi '_x)^2-{\varphi '_x}^2)+i(\varphi '(x)(hD+\psi ')+(hD+\psi 
')\varphi '))\big) a\nonumber\\
&=&(hL-h^2\Delta )a,\nonumber
\end{eqnarray} 
where $L$ is the transport \op{} given by
\ekv{sp.15}
{
L=\psi 'D+D\psi '+i(\varphi 'D+D\varphi ').
}
Along the projection of each bicharacteristic leaf this is an elliptic 
\op{} of Cauchy-Riemann type and if we assume that the leaves are open and 
simply connected, then (see \cite{DuHo}) there exists a non-vanishing 
smooth \fu{} $a\in C^\infty $ such that 
\ekv{sp.16}
{
La=0.}

\par Recall that $q\in L^\infty (\widetilde{\Omega })$. Assume that $a$ 
in (\ref{sp.16}) is well-defined in a \neigh{} of $\overline{\Omega }$. 
Then from (\ref{sp.14}), we see that with $P=P_0+h^2q$:
\ekv{sp.17}
{
P e^{{1\over h}(-\varphi +i\psi )}a=e^{-\varphi /h}h^2d,
}
with $d={\cal O}(1)$ in $L^\infty $ and hence in $L^2$. Now apply 
Proposition \ref{PropSP.1} with $\varphi $ replaced by $-\varphi $, to find 
$r\in H^1(\Omega )$ with $h\Vert r\Vert _{H^1}\le Ch^2$, such that 
$$e^{\varphi /h}Pe^{-\varphi /h}e^{i\psi /h}r=-h^2d,$$
i.e.  \ekv{sp.18}
{
P(e^{{1\over h}(-\varphi +i\psi )}(a+r))=0.
}

\section{More use of the Carleman estimate}\label{SectionUse}
\setcounter{equation}{0}

In Section \ref{SectionCe} we derived a Carleman estimate for $e^{\varphi 
/h}u$ 
when $h^2(-\Delta +q)u=v$ when $\varphi $ is a smooth limiting Carleman 
weight 
with non-vanishing gradient. In order to stick close to the paper
\cite{BuUh}, we write the corresponding estimate for $e^{-\varphi 
/h}u$, when $(-\Delta +q)u=v$, ${u_\vert}_{\partial \Omega }=0$:
\begin{eqnarray}\label{use.1}
&
{h^3\over C_0}((\varphi '_x\cdot \nu )e^{-\varphi /h}\partial _\nu u\vert 
e^{-\varphi /h}\partial _\nu u)_{\partial \Omega _+}+{h^2\over C_0}(\Vert 
e^{-\varphi /h}u\Vert ^2+\Vert e^{-\varphi /h}h\nabla u\Vert ^2)&\\
&\le h^4\Vert e^{-\varphi /h}v\Vert ^2-C_0h^3((\varphi '_x\cdot 
\nu )e^{-\varphi /h}\partial _\nu u\vert e^{-\varphi /h}\partial _\nu 
u)_{\partial \Omega _-},
&\nonumber
\end{eqnarray}
where $\nu $ is the exterior unit normal and $\Omega _{\pm}=\{ x\in 
\partial \Omega ;\, \pm \nu \cdot \varphi '>0\}$.

\par Let $q_1,q_2\in L^\infty (\Omega )$ be two potentials. Let 
\ekv{use.2}
{
u_2=e^{{1\over h}(\varphi +i\psi _2)}(a_2+r_2(x;h)),\hbox{ with }(\Delta 
-q_2)u_2=0,\ \Vert r_2\Vert _{H^1}={\cal O}(h).
}
Here $\psi _2$ is chosen as in Section \ref{SectionSP} so that 
$(\varphi ')^2=(\psi _2')^2=\varphi '\cdot \psi 
_2'=0$ and so that the integral leaves of the commuting 
\vf{}s $\varphi '\cdot \partial _x,\psi 
'\cdot \partial _x$ are simply connected in $\Omega $. $a_2$
is smooth in a \neigh{} of $\overline{\Omega }$ and everywhere 
non-vanishing.

\par Let ${\cal N}_q$ be the Dirichlet to Neumann map for the potential 
$q$ and let 
\begin{eqnarray*}
\partial \Omega _{-,\epsilon _0}&=&\{ x\in\partial \Omega ;\, \nu(x) \cdot 
\varphi 
'_x(x)<\epsilon _0\},\\
\partial \Omega _{+,\epsilon _0}&=&\{ x\in \partial \Omega ;\, \nu(x) 
\cdot 
\varphi '_x(x)\ge \epsilon _0\} ,
\end{eqnarray*}
for some fixed $\epsilon _0>0$, so that $\partial \Omega _{+,\epsilon 
_0}\subset \partial \Omega _+$, $\partial \Omega _-\subset \partial \Omega 
_{-,\epsilon _0}$. Here $\nu(x)$ denotes the unit outer normal to $\partial
\Omega.$

\par Assume
\ekv{use.3}
{{\cal N}_{q_1}(f)={\cal N}_{q_2}(f),\hbox{ in }\partial \Omega 
_{-,\epsilon _0},\hbox{ for all }f\in H^{1\over 2}(\partial \Omega ).}

\par Let $u_1\in H^1(\Omega )$ solve 
\ekv{use.3.5}
{(\Delta -q_1)u_1=0,\ {{u_1}_\vert}_{\partial \Omega }=
{{u_2}_\vert}_{\partial \Omega }.}
Then by the assumption (\ref{use.3}), we have \ekv{use.4}{\partial _{\nu 
}u_1=
\partial _\nu u_2\hbox{ in }\partial \Omega 
_{-,\epsilon _0}.
}
Put $u=u_1-u_2$, $q=q_2-q_1$, so that 
\ekv{use.5}
{
{\rm supp\,}({\partial _\nu u_\vert }_{\partial \Omega })\subset \partial 
\Omega _{+,\epsilon _0},
}
and 
\ekv{use.6}
{
(\Delta -q_1)u=(\Delta -q_1)u_2=qu_2,\ {u_\vert}_{\partial \Omega }=0.
}
For $v\in H^1(\Omega )$ with $\Delta v\in L^2(\Omega )$, we get using 
(\ref{use.5}),(\ref{use.6}) and Green's formula:
\ekv{use.7}
{
\int_\Omega qu_2\overline{v}dx=\int_\Omega (\Delta 
-q_1)u\overline{v}dx=\int_\Omega u\overline{(\Delta 
-\overline{q_1})v}dx+\int_{\partial \Omega _{+,\epsilon _0}}(\partial _\nu 
u)
\overline{v}\,S(dx).}

\par As in Section \ref{SectionSP} we can construct
\ekv{use.8}
{v=e^{-{1\over h}(\varphi +i\psi _1)}(a_1+r_1),}
 with $\psi _1$ satisfying $\varphi '\cdot \psi _1'=0$, $(\varphi ')^2=(\psi 
 _1')^2$, with $a_1(x)$ non-vanishing and smooth, 
and with $\Vert r_1\Vert _{H^1(\Omega )}={\cal O}(h)$, so that 
\ekv{use.9}
{(\Delta -\overline{q}_1)v=0.}
Then (\ref{use.7}) becomes
\ekv{use.10}
{
\int_\Omega q e^{{i\over h}(\psi _1+\psi 
_2)}(a_2+r_2)\overline{(a_1+r_1)}dx=\int_{\partial \Omega _{+,\epsilon 
_0}}(\partial _\nu u)e^{-{1\over h}(\varphi -i\psi _1)}\overline{(a_1+r_1)}\, 
S(dx).
}
We shall work with $\psi_1, \psi _2, \varphi $
slightly $h$-dependent in such a way that 
\ekv{use.10.5}
{{1\over h}(\psi _1+\psi _2)\to f,\ h\to 0.}
Recall that 
\ekv{use.11}
{\Vert r_j\Vert_{H^1}={\cal O}(h).}
Then using that $q\in L^\infty $, we see that the \lhs{} of (\ref{use.10}) 
converges to
\ekv{use.12}
{\int_\Omega a_2\overline{a}_1q(x)e^{if(x)}dx.}

\par For the \rhs{} of (\ref{use.10}), we have, using (\ref{use.1}), 
for $(\Delta -q_1)$ and (\ref{use.6}):
\begin{eqnarray}\label{use.13}
&& \vert \int_{\partial \Omega _{+,\epsilon _0}}(\partial _\nu u) 
e^{-{1\over h}(\varphi -i\psi _1)}(\overline{a_1+r_1})\,S(dx)\vert ^2\\
&\le & \Vert a_1+r_1\Vert ^2_{\partial \Omega _{+,\epsilon _0}} 
\int_{\partial \Omega _{+,\epsilon _0}}(e^{-\varphi /h}\vert \partial _\nu u
\vert )^2S(dx)\nonumber\\
&\le & \Vert a_1+r_1\Vert ^2_{\partial \Omega _{+,\epsilon _0}} 
{1\over \epsilon _0}\int_{\partial \Omega _{+,\epsilon _0}}(\varphi '\cdot 
\nu )(e^{-\varphi /h}\vert \partial _\nu u\vert )^2 S(dx)\nonumber\\
&\le & {1\over \epsilon _0}\Vert a_1+r_1\Vert _{\partial \Omega 
_{+,\epsilon _0}}^2(C_0h\Vert e^{-\varphi /h}qu_2\Vert 
^2-C_0^2\int_{\partial \Omega _-}(\varphi '\cdot \nu )(e^{-\varphi 
/h}\vert \partial _\nu u\vert )^2S(dx)).\nonumber
\end{eqnarray}
Here $\partial _\nu u=0$ on $\partial \Omega _-$, and using also 
(\ref{use.2}), we get
\ekv{use.14}
{
\vert \int_{\partial \Omega _{+,\epsilon _0}}(\partial _\nu u)e^{-{1\over 
h}(\varphi -i\psi _1)}(\overline{a_1+r_1})S(dx)\vert ^2\le {C_0h\over 
\epsilon _0}
\Vert a_1+r_1\Vert ^2_{\partial \Omega 
_{+,\epsilon _0}}\Vert q(a_2+r_2)\Vert ^2.}

\par Here $\Vert q(a_2+r_2)\Vert ^2={\cal O}(1)$, by (\ref{use.11}). 
Since $r_1={\cal O}(h)$ in the semiclassical $H^1$-norm, we have 
$r_1={\cal 
O}(1)$ in the standard ($h=1$) $H^1$-norm. Hence 
\ekv{use.15}
{
{{r_1}_\vert}_{\partial \Omega }={\cal O}(1)\hbox{ in }L^2.
}
Consequently, the \rhs{} of (\ref{use.10}) tends to 0, when $h\to 0$, 
and letting $h\to 0$ there, we get 
\ekv{use.16}
{
\int_{\Omega }q(x)a_2(x)\overline{a}_1(x) e^{if(x)}dx=0,
}
for all $f$ that can be attained as limits in (\ref{use.10.5}).

\par Finally, we remark that if $\varphi $ is real-\an{}, then in the above 
constructions, we may arrange so that $\psi _j$
and $a_j$ have the same property.

\section{End of the proof of Theorem 1.2}\label{SectionPf}
\setcounter{equation}{0}

\par From now on, we assume that the dimension $n$ is $\ge 3$. 
We choose $\varphi (x)=\ln \vert x-x_0\vert $ for $x_0$ varying in a 
small open set separated from $\overline{\Omega }$ by some fixed affine 
hyperplane $H$. Notice that $\varphi$ is a limiting Carleman
weight in the sense of (\ref{ce.4}). 
We need a \sufly{} rich \fy{} of \fu{}s $f$ in 
(\ref{use.16}) and recall that these \fu{}s are the ones that appear in 
(\ref{use.10.5}) with $\psi _j$ analytic near $\overline{\Omega }$ and 
satisfying $(\psi _j')^2=(\varphi ')^2$, $\psi _j'\cdot \varphi '=0$. Changing 
the sign of $\psi _2$ we can also view $f$ as a limit ${1\over h}(\psi 
_1-\psi _2)$ for suitable such $h$-dependent \fu{}s $\psi _j$. More 
precisely, we can take an analytic \fy{} $\psi (x,\alpha )$ depending on 
the additional parameters $\alpha =(\alpha _1,...,\alpha _k)$, with 
$\psi (\cdot ,\alpha )$ \sy{}ing 
\ekv{pf.1}
{(\psi _x')^2=(\varphi '_x)^2,\ \psi _x'\cdot \varphi _x'=0,}
and then take
\ekv{pf.2}
{f(x)=\langle \psi '_\alpha (x,\alpha ),\nu (\alpha )\rangle ,}
where $\nu (\alpha )$ is a tangent vector in the $\alpha $-variables.

\par We first discuss the choice of $\psi $. Since $\varphi '_x$ is radial, 
\wrt{} $x_0$, the second condition in (\ref{pf.1}) means that $\psi (x)$ 
is 
positively \hg{} of degree $0$ \wrt{} $x-x_0$. A necessary and sufficient 
condition for $\psi $ (at least if we work in some cone with vertex at 
$x_0$) is then that
\ekv{pf.3}
{
(\psi _x')^2=(\varphi '_x)^2,}
on a suitable open subset $x_0+r_0W$ of $x_0+r_0S^{n-1}$, for some fixed 
$r_0>0$. The necessity is obvious and the sufficiency follows easily by 
extending $\psi $ to be a positively \hg{} \fu{} of degree $0$ in the 
variables $x-x_0$. 

\par Here is an explicit choice of a suitable open set in (\ref{pf.3}): 
Let $r_0>0$ be large enough so that $\overline{\Omega }\subset 
B(x_0,r_0)$. 
Let $x_0+r_0W\subset \partial B(x_0,r_0)$ be defined by 
\ekv{pf.4}
{x_0+r_0W=\partial B(x_0,r_0)\cap H_+,}
where $H_+$ is the open half-space delimited by the affine hyper-plane 
$H$, 
for which $x_0\notin H_+$ (so that $\overline{\Omega }\subset H_+$). Then 
$\overline{\Omega }$ is contained in the open cone $x_0+{\bf R}_+W$, so if 
we choose $\psi $ on 
$x_0+r_0W$ as in (\ref{pf.3}) and 
extend by homogeneity, we know that $\psi $ will be smooth near 
$\overline{\Omega }$.

\par Let $y_0\in \partial B(0,1)\setminus\overline{W}$ be 
such that the
antipodal point $-y_0$ also is outside $\overline{W}$ and  
define
\ekv{pf.5}
{\psi (x,y)=d_{S^{n-1}}(x,y).}
Then $\psi \in C^\infty (W\times {\rm neigh\,}(y_0))$ and the function 
$\psi ((x-x_0)/\vert x-x_0\vert ,y)\in C^\infty (\overline{\Omega 
}\times {\rm neigh\,}(y_0))$ will  \sy{} (\ref{pf.1}). Since the domain of 
definition does not contain antipodal points, we remark that 
\ekv{pf.6}
{
\psi _{x,y}''\hbox{ is of rank }n-2\hbox{ and }{\cal R}(\psi 
''_{x,y})=(\psi '_x)^\perp , {\cal N}(\psi ''_{x,y})=(\psi '_y).
}
This follows from basic properties of the geodesic flow (and remains true 
more generally for $\psi (x,y)=d(x,y)$ on a Riemannian \mfld{} as long 
as $x,y$ are not conjugate points.)

\par For $x\in W\subset S^{n-1}$, $(y,\nu )\in TS^{n-1}$, $y\in{\rm 
neigh\,}(y_0)$, we put
\ekv{pf.7}
{\widetilde{f}(x;y,\nu )=\psi '_y(x,y)\cdot \nu .}
Then 
\ekv{pf.8}
{
\widetilde{f}'_x(x;y,\nu )=\psi ''_{x,y}(x,y)(\nu ).
}
In view of (\ref{pf.6}), we see that this vanishes precisely when $\nu 
\parallel \psi '_y(x,y)$, i.e. when $\nu $ is parallel to the (arrival) 
direction of the minimal geodesic from $x$ to $y$. Restricting $\nu $ to 
non-vanishing directions which are close to be parallel to the plane $H$, 
we can assure that 
\ekv{pf.8.5}
{\widetilde{f}'_x(x;y,\nu )\ne 0.}

\begin{lemma}\label{LemmaPf.1}
$\widetilde{f}''_{x,(y,\nu )}$ has maximal rank $n-1$.
\end{lemma}

\begin{proof} We already know that $\widetilde{f}''_{x,\nu }=\psi 
''_{x,y}$ is of rank $n-2$ and that the image of this matrix is equal to 
$(\psi '_x)^\perp$. Consequently, we consider 
$$g(y)=\psi _x'(x,y_0)\cdot \psi ''_{x,y}(x,y)(\nu )=\langle \psi 
''_{x,y}(x,y)\vert \psi '_x(x,y_0)\otimes \nu \rangle $$
as a \fu{} of $y\in{\rm neigh\,}(y_0)$. The function vanishes for $y=y_0$ 
and can also be written 
$$\langle \psi ''_{x,y}(x,y)\vert (\psi '_x(x,y_0)-\psi '_x(x,y))\otimes 
\nu \rangle =\langle \psi ''_{x,y}(x,y)(\nu )\vert \psi 
''_{x,y}(x,y)(y_0-y)\rangle +{\cal O}((y_0-y)^2).$$
From this expression, we see that the $y$-differential is non-vanishing  
and 
hence the range of 
$\widetilde{f}''_{x,(y,\nu )}$ contains vectors that are not 
orthogonal to $\psi '_x(x,y)$.
\end{proof}

\par Now consider 
\ekv{pf.9}
{\Psi (x;y,\widetilde{x})=\psi ({x-\widetilde{x}\over \vert 
x-\widetilde{x}\vert },y)\in C^\infty (\overline{\Omega }
\times {\rm neigh\,}(y_0,S^{n-1}) \times {\rm 
neigh\,}(x_0,{\bf R}^n)).
}
$\Psi $ is analytic, real and satisfies (\ref{pf.1}) with $\varphi 
(x)=\Phi (x,\widetilde{x})=\ln \vert x-\widetilde{x}\vert $. 
We can take $\alpha =y$
and (\ref{pf.2}) becomes 
\ekv{pf.10}
{f(x)=f(x;\theta )=\langle \Psi '_y | \nu \rangle ,\ 
\theta =(y,\widetilde{x},\nu),}
with $(y,\nu )\in TS^{n-1}$. Lemma \ref{LemmaPf.1} shows that 
$f''_{x,(y,\nu )}$ has rank $n-1$ and indeed the image of this matrix is 
the tangent space of $\partial B(\widetilde{x},\vert x-\widetilde{x}\vert 
)$ 
at $x$. 
Since $f'_x$ is a non-vanishing element of $T_x(\partial B(\widetilde{x},
|x-\widetilde{x}|))$, we can vary $\widetilde{x}$ infinitesimally to see 
that $f''_{x,\widetilde{x}}(\widetilde{\mu })\not\in T_x(\partial 
B(\widetilde{x},|x-\widetilde{x}|))$ for a suitable $\widetilde{\mu } 
\in {\bf R}^n$. It is 
then clear that 
\ekv{pf.11}
{f''_{x,\theta }=
f''_{x,(y,\widetilde{x},\nu )}\hbox{ has maximal rank }n,}
and hence that the map
\ekv{pf.12}
{
{\rm neigh\,}(\overline{\Omega })\ni x\mapsto f'_\theta (x,\theta )\in 
{\bf R}^{3n-2}
}
has injective differential. 

\begin{lemma}\label{LemmaPf.2}
The map (\ref{pf.12}) is injective.
\end{lemma}

\begin{proof}
Let $x_1,x_2\in{\rm neigh\,}(\overline{\Omega })$ be two
points with 
\ekv{pf.13}
{
f'_\theta (x_1,\theta )=f'_\theta (x_2,\theta ),
}
for some $\theta =(y,\widetilde{x},\nu )$. 
Taking the $\nu 
$-component of this relation, we get
\ekv{pf.14}
{
\psi '_y(\widetilde{x}_1,y)=\psi '_y(\widetilde{x}_2,y),\ 
\widetilde{x}_j={x_j-\widetilde{x}\over \vert x_j-\widetilde{x}\vert }.}

\par This means that $\widetilde{x}_1$, $\widetilde{x}_2$, $y$ belong 
to the same geodesic $\gamma $ and this geodesic is minimal (i.e. distance 
minimizing) on some segment that contains these three points in its 
interior. If $\widetilde{x}_1\ne \widetilde{x}_2$, we may assume that 
$d(\widetilde{x}_2,y)<d(\widetilde{x}_1,y)$. For $y\in{\rm 
neigh\,}(y_0,S^{n-1})$, we have 
$$
d(\widetilde{x}_1,\widetilde{x}_2)+d(\widetilde{x}_2,y)- 
d(\widetilde{x}_1,y)=:g(y), \ g(y)\sim d(y,\gamma )^2.
$$
It follows that
$$f(\widetilde{x}_2;y,\widetilde{x},\nu )-
f(\widetilde{x}_1;y,\widetilde{x},\nu )=g'(y)\cdot \nu ,$$
and using that $\nu $ is not parallel to $\dot{\gamma }$ at $y_0$, we see 
that this \fu{} has a non-vanishing $y$-gradient at $y=y_0$, in 
contradiction with (\ref{pf.13}). Thus, $\widetilde{x}_1=\widetilde{x}_2$, 
or in other words, $x_1$ and $x_2$ belong to the same half-ray through 
$\widetilde{x}$.

\par Taking the $\widetilde{x}$-component of (\ref{pf.13}), we get
$${{\nabla _{\widetilde{x}}\langle \psi '_y ({x-\widetilde{x}\over \vert 
x-\widetilde{x}\vert },
y ),\nu \rangle }_\vert}_{x=x_1}={{\nabla _{\widetilde{x}}\langle 
\psi '_y ({x-\widetilde{x}\over \vert 
x-\widetilde{x}\vert },
y ),\nu \rangle }_\vert}_{x=x_2}
.$$
These quantities are clearly non-vanishing and if $x_1\ne x_2$, they 
differ by a factor $\ne 1$, since $\nabla 
_{\widetilde{x}}({x-\widetilde{x}\over \vert 
x-\widetilde{x}\vert })$ is \hg{} of degree $-1$ in $x-\widetilde{x}$. 
Thus $x_1=x_2$. 
\end{proof}

\par Now apply (\ref{use.16}) with $f(x)=f(x,\theta )$:
\ekv{pf.15}
{
\int e^{if(x,\theta )}a_2\overline{a}_1q(x)dx =0,
}
where $a_2$, $a_1$ are analytic non-vanishing functions of
$x,y,\widetilde{x}$ in a neighborhood of $\overline{\Omega }\times \{
y_0\} \times \{ x_0\} $.
Since $f(x,\theta )=f(x;y,\widetilde{x},\nu )$ depends linearly on 
$\nu $, we can replace $\nu $ by $\lambda \nu $ and 
get 
\ekv{pf.16}
{
\int e^{i\lambda f(x,\theta )}a_2\overline{a}_1q(x)dx=0,\
\lambda \ge 1.
}
Now represent $\theta $ by some analytic real coordinates
$\theta 
_1,\theta _2,...,\theta _N$ near some fixed given point $\theta
_0=(y_0,x_0,\nu _0)$. If $x,z\in\overline{\Omega }$, $w\in{\rm
neigh\,}(\theta _0)$, we consider the \fu{}
\ekv{pf.17}
{
\theta \mapsto -f(z,\theta )+f(x,\theta )+{i\over 2}(\theta
-w)^2.
}
For $x=z$, we have the unique \nondeg{} critical point $\theta
=w$, 
while for $x\ne z$ there is no real critical point in view of Lemma
\ref{LemmaPf.2}. For $x\approx z$ we 
have a unique complex critical
point which is close to $w$, and we 
introduce the corresponding critical
value 
\ekv{pf.17.5}
{\psi (z,x,w)={\rm v.c.}_\theta  (-f(z,\theta )+f(x,\theta
)+{i\over 
2}(\theta -w)^2).}
From (\ref{pf.12}) and standard 
estimates
on critical values in connection with the complex stationary 
phase method
(\cite{MeSj, Sj}), we deduce that 
\ekv{pf.18}
{
\Im \psi (z,x,w)\sim (z-x)^2,\ z,x\in\overline{\Omega },\,\, z\approx x.
}  
Moreover, when $x=z$, we have 
\ekv{pf.19}
{
\psi '_z(z,z,w)=-f'_z(z,w),\ \psi '_x(z,z,w)=f'_z(z,w),\ \psi (z,z,w)=0.}

\par We now multiply (\ref{pf.16}) by $\chi (\theta -w)e^{i\lambda {i\over
2}(\theta -w)^2-i\lambda f(z,\theta )}$, and integrate \wrt{} $\theta $,
to get \ekv{pf.20}
{
\int e^{i\lambda \psi (z,x,w)}a(z,x,w;\lambda )\chi (z-x)q(x)dx={\cal 
O}(e^{-{\lambda \over C}}).
} 
Here $\chi $
denote (different) standard cutoffs to a \neigh{} of $0$, and $a$ is an 
elliptic classical analytic symbol of order 0.

\par Now restrict $w$ to an $n$-dimensional \mfld{} $\Sigma $ which passes
through $\theta _0$, and write $(z,-f'_z(z,\theta ))=(\alpha _x,\alpha
_\xi )=\alpha $. Then we rewrite (\ref{pf.20}) as
\ekv{pf.21}
{
\int e^{{i\lambda }\psi (\alpha ,x)}a(\alpha ,x;h)\chi (\alpha 
_x-x)q(x)dx={\cal O}(e^{-{\lambda \over C}}),
}
implying that
\ekv{pf.22}
{
(z,-f'_z(z,\theta _0))\notin {\rm WF}_a(q),
} 
since we can apply the standard FBI-approach (\cite{Sj}). Notice that 
(\ref{pf.18}), (\ref{pf.19}) give: 
\ekv{pf.23}
{
\psi (\alpha ,x)=(\alpha _x-x)\cdot \alpha _\xi +{\cal O}((\alpha 
_x-x)^2),\ \Im \psi (\alpha ,x)\sim (\alpha _x-x)^2,}
and we can choose $\Sigma $ so that the map ${\rm neigh\,}(z_0)\times
\Sigma \ni (z,\theta )\mapsto (z,-f'_z(z,\theta ))$ is local
diffeomorphism near any given fixed point $z_0\in \overline{\Omega }$.
\medskip

\par\noindent 
\it End of the proof of \Th{} \ref{Th0.2}. \rm Fix $\theta _0$ as above, 
so 
that $0\ne -f'_z(z,\theta _0)\notin {\rm WF}_a(q)$ for all $z$ in some 
\neigh{} of $\overline{\Omega }$. (Notice that $q$ now denotes the 
extension by $0$ of the originally defined function on $\Omega $.) Let 
$z_0$ be a point in ${\rm supp\,}(q)$, where ${f(\cdot ,\theta _0)_\vert 
}_{{\rm supp\,}(q)}$ is minimal. Then $-f'_z(z_0,\theta _0)$ belongs to 
the 
exterior conormal cone of ${\rm supp\,}(q)$ at $z_0$ and we get a 
contradiction between (\ref{pf.22}) and the fact that all such exterior 
conormal directions have to belong to $WF_a(q)$. (This is the
wavefront version of Holmgren's uniqueness theorem, 
due to H{\"o}rmander (\cite{Ho}) and 
Sato-Kawai-Kashiwara (remark by Kashiwara in \cite{SaKaKa}).) 
\hfill$\Box$\medskip 

\section{Complex geometrical optics solutions with Dirichlet data on part 
of
the boundary}\label{SectionCo}
\setcounter{equation}{0}

In this section we prove Theorem \ref{Th0.1}.

We first use the Carleman estimate (\ref{ce.24.7}) and the Hahn-Banach 
theorem
to construct CGO solutions for the conjugate operator
$P^\ast_\varphi=(e^{\frac{\varphi}{h}}Pe^{-\frac{\varphi}{h}})^*$ where 
$*$
denotes the adjoint. Notice that $P^*_\varphi$ has the same form as
$P_\varphi$ except that $q$ is replaced by $\bar q$ and $\varphi $ by $-\varphi $.

\begin{prop}\label{PropCo.1}
Let $\varphi$ be as in (\ref{ce.24.7}).
Let $v\in H^{-1}(\Omega)$, $v_-\in L^2(\partial\Omega_-; (-\varphi'\cdot 
\nu)S(dx))$. Then
$\exists\,u\in H^0(\Omega)$ such that $$P^\ast_\varphi u=v, \quad
u\bigr|_{\partial\Omega_-}=v_-.$$
Moreover
\begin{equation}\label{co.1}
\|u\|_{H^0}+\sqrt{h}\|(\varphi'\cdot 
\nu)^{-\frac{1}{2}}u\|_{\partial\Omega_+}
\le C(\frac{1}{h}\|v\|_{H^{-1}}+\sqrt{h}\Vert(-\varphi'\cdot
\nu)^{-\frac{1}{2}}v_-\Vert_{\partial\Omega_-}).
\end{equation}
\end{prop}

\begin{proof}
We use the Carleman estimate (\ref{ce.24.7}). Let $v$ as in 
the proposition. For
$w\in(H^1_0\cap H^2)(\Omega)$ we have $$|(w|v)_\Omega + (h\partial_\nu
w|v_-)_{\partial\Omega_-}| \le \|w\|_{H^1} \|v\|_{H^{-1}}+\left 
( (-\varphi'\cdot
\nu)^{\frac{1}{2}}h\partial_\nu w|(-\varphi'\cdot 
\nu)^{-\frac{1}{2}}v_-\right )_{\partial \Omega _-}.$$
Therefore
$$|(w,v)_\Omega+(h\partial_\nu w|v_-)_{\partial\Omega_-}|\le C
(\frac{1}{h}\|v\|_{H^{-1}}h\| w\|_{H^{1}}+\frac{1}{\sqrt{h}}\| (-\varphi'\cdot
\nu)^{-\frac{1}{2}}v_-\|_{\partial\Omega_-}\sqrt{h}\|
(-\varphi'\cdot \nu)^{\frac{1}{2}}h\partial_\nu w\|_{\partial\Omega_-}).$$ 
Now by using (\ref{ce.24.7}) we get  $$|(w,v)_\Omega + (h\partial_\nu
w|v_-)_{\partial\Omega_-}|\le
C(\frac{1}{h}\|v\|_{H^{-1}}+\frac{1}{\sqrt{h}}\|(-\varphi'\cdot
\nu)^{-\frac{1}{2}}v_-\|_{\partial\Omega_-})(\| P_\varphi
w\|+\sqrt{h}\|(\varphi'\cdot \nu)^{\frac{1}{2}}h\partial_\nu 
w\|_{\partial\Omega_+}).$$ By the Hahn-Banach theorem, $\exists\, u\in
H^0(\Omega)$, $u_+\in L^2(\partial\Omega_+, (\varphi'\cdot
\nu)^{-\frac{1}{2}}dS)$, $u_+$ on $\partial\Omega_+$ such that
\begin{equation}\label{co.2}
(w,v)_\Omega + (h\partial_\nu w|v_-)_{\partial\Omega_-} = (P_\varphi
w|u)+(h\partial_\nu w|u_+)_{\partial\Omega_+},\quad \forall\,w\in(H^1_0\cap
H^2)(\Omega)
\end{equation}
with
\begin{equation}\label{co.3}
\|u\|_{H^0}+\frac{1}{\sqrt{h}}\|(\varphi'\cdot
\nu)^{-\frac{1}{2}}u_+\|_{\partial\Omega_+} \le
C(\frac{1}{h}\|v\|_{H^{-1}}+\frac{1}{\sqrt{h}}\|(-\varphi'\cdot
\nu)^{-\frac{1}{2}}v_-\|_{\partial\Omega_-}).
\end{equation}
Since $P_\varphi=-h^2\Delta+$ a first order operator, and
$w\bigr|_{\partial\Omega}=0$ we have
$(P_\varphi w|u)=(w|P^*_\varphi u)-h^2(\partial_\nu
w|u)_{\partial\Omega}$. 

Using this in (\ref{co.2}) we obtain
$$0=(w|v-P^*_\varphi u) +  h((\partial_\nu
w|1_{\partial\Omega_-}v_-)_{\partial\Omega} - (\partial_\nu
w|1_{\partial\Omega_+} u_+)_{\partial\Omega} + (\partial_\nu
w|hu)_{\partial\Omega})$$ where $1_{\partial\Omega_\pm}$ denotes the 
indicator
function of $\partial\Omega_\pm$.

By varying $w$ in $(H^1_0\cap H^2)(\Omega)$ we get
$$P^*_\varphi u=v, \quad
hu\bigr|_{\partial\Omega}=-1_{\partial\Omega_-}v_-+1_{\partial\Omega_+}u_+.$$
which implies the proposition after replacing $v_-$ above by $-hv_-$.
\end{proof}

Let $W_-\subset\partial \Omega_-$ be an arbitrary strict open subset of
$\partial\Omega_-$. We next want to modify the choice of $u_2$ in
(\ref{use.2}) so that $u_2\bigr|_{W_-}=0$.

\begin{prop}\label{PropCo.2}
Let $a_2$, $\varphi$, $\psi_2$ be as in (\ref{use.2}). Then we can 
construct
a solution of
\begin{equation}\label{co.4}
P\widetilde u_2=0, \quad \widetilde u_2\bigr|_{\overline W_-}=0
\end{equation}
of the form 
\begin{equation}\label{co.5}
\widetilde u_2=e^{\frac{1}{h}(\varphi+i\psi_2)}(a_2+\widetilde r_2)+u_r
\end{equation}
where $u_r=e^{i\frac{l}{h}}b(x;h)$ with $b$ a symbol of order
zero in $h$ and
\begin{equation}\label{co.6}
\Im l(x)=-\varphi (x)+k(x)
\end{equation}
where $k(x)\sim\dist(x,\partial\Omega_-)$ in a neighborhood of
$\partial\Omega_-$ and $b$ has its support in that \neigh{}. Moreover, 
$\Vert \widetilde{r}_2\Vert _{H^0}={\cal O}(h)$, 
$\widetilde{r}_2\bigr|_{\partial \Omega _-}=0$, $\Vert (\varphi '\cdot 
\nu )^{-{1\over 2}}\widetilde{r}_2\Vert _{\partial \Omega _+}={\cal 
O}(h^{1\over 2})$.
\end{prop}

\begin{proof}
We start by constructing a WKB solution $u$ in $\Omega$ of
\begin{equation}\label{co.7}
-h^2\Delta u=0,
u\bigr|_{\partial\Omega_-}=e^{\frac{1}{h}(\varphi+i\psi_2)}(\chi a_2)
\bigr|_{\partial\Omega_-}
\end{equation}
where $\chi\in C^\infty_0(\partial\Omega_-)$, $\chi=1$ on $\bar W_-$.

We try
$u=e^{\frac{i}{h}l(x)}b(x;h)$. The eikonal equation for $l$ is
\begin{equation}\label{co.8}
\begin{array}{rcl}
(l')^2 & =  &0 \mbox{ to infinite order at }\partial\Omega \\
\left.l\right|_{\partial\Omega_-} & = & \psi_2-i\varphi.
\end{array}
\end{equation}

Of course $g:=\psi_2-i\varphi$ is a solution but we look for the second
solution, corresponding to having $u$ equal to a ``reflected wave''. We 
decompose on $\partial\Omega_-$
\[
g'=g'_t+g'_\nu
\]
where $t$ denotes the tangential part and $\nu$ the normal part.

Then in order to satisfy the eikonal equation we need 
\[
0=(g'_t)^2+(g'_\nu)^2.
\]
Therefore we can solve (\ref{co.8}) to $\infty$-order at $\partial\Omega_-$
with $l$ satisfying
\[
l\bigr|_{\partial\Omega_-} = g\bigr|_{\partial\Omega_-}, \partial_\nu
l\bigr|_{\partial\Omega_-} = - \partial_\nu g\bigr|_{\partial\Omega_-}.
\]
By the definition of $\partial\Omega_-$ we have
\[
\partial_\nu \Im g=-\partial_\nu \varphi>0\mbox{ on }\partial\Omega_-.
\]
Since $\nu$ is the unit exterior normal we have that (\ref{co.6}) is
satisfied.

Solving also the transport  equation to $\infty$-order, at the boundary we
get a symbol $b$ of order $0$ with support \ably{} close to $\supp \chi $, such that
\[
\left\{\begin{array}{rcl}
-h^2\Delta(e^{\frac{il}{h}}b(x;h)) & = &
e^{\frac{il}{h}}{\cal O}((\dist(x,\partial\Omega))^\infty+h^\infty) \\
\left.e^{\frac{il}{h}}\right|_{\partial\Omega} & = &
\left.e^{\frac{ig}{h}}\chi a_2\right|_{\partial\Omega}.
\end{array}\right.
\]
Our new WKB input to $u_2$ will be
\[
(e^{\frac{ig}{h}}a_2-e^{i\frac{l}{h}}b).
\]
Instead of (\ref{sp.17}) we get
\begin{equation}\label{co.9}
P(e^{i\frac{g}{h}}a_2-e^{\frac{il}{h}}b) = e^{\frac{\varphi}{h}}h^2d
\end{equation}
where $d={\cal O}(1)$ in $L^2(\Omega)$.

Using Proposition \ref{PropCo.1} we can solve
\begin{eqnarray*}
e^{-\frac{\varphi}{h}}Pe^{\frac{\varphi}{h}}(e^{i\frac{\psi_2}{h}}\widetilde
r_2)
& = & -h^2d \\
\left.\widetilde r_2\right|_{\partial\Omega_-} & = & 0
\end{eqnarray*}
with
\[
\|\widetilde r_2\|_{H^0} + \sqrt{h}\|(\varphi'\cdot
\nu)^{-\frac{1}{2}}\widetilde r_2\|_{\partial\Omega_+} \le 
\frac{C}{h}\|h^2d\|_{H^{-1}}
= {\cal O}(h).
\]
Thus
\begin{equation}\label{co.10}
\|\widetilde r_2\|={\cal O}(h), \quad \|(\varphi'\cdot
\nu)^{-\frac{1}{2}}\widetilde r_2\|_{\partial\Omega_+}={\cal O}(\sqrt{h}).
\end{equation}
Now we take
\begin{equation}\label{co.11}
u_2=e^{\frac{1}{h}(\varphi+i\psi_2)}(a_2+\widetilde r_2)-e^{\frac{il}{h}}b.
\end{equation}
Clearly $Pu_2=0$, $u_2\bigr|_{\partial\Omega}=0$ in $\overline W_-$.
\end{proof}

\emph{Proof of Theorem \ref{Th0.1}}. Let $u_2=\widetilde{u}_2$ be as in Prop \ref{PropCo.2}. Let
$u_1\in H^1(\Omega)$ solve (\ref{use.3.5})
\[
(\Delta-q_1)u_1=0, \quad u_1\bigr|_{\partial\Omega} =
u_2\bigr|_{\partial\Omega}.
\]
By construction we have that $\supp\, u_i\bigr|_{\partial \Omega }\cap\overline{W_-}=\emptyset$,
$i=1,2$. As in Section \ref{SectionUse}, let $u=u_1-u_2$, $q=q_1-q_2$.
Then  (\ref{use.5}) and (\ref{use.6}) are valid and in fact $u\in
H^2(\Omega)$ so that the Green's formula (\ref{use.7}) is also valid. Now
choose $v$ as in (\ref{use.8}), (\ref{use.9}). Then instead of (\ref{use.10}) 
we get
\begin{equation}\label{co.12}
\begin{array}{c}
\int_\Omega 
qe^{\frac{i}{h}(\psi_1+\psi_2)}\overline{(a_1+r_1)}(a_2+\widetilde r_2)dx -
\int_\Omega qe^{\frac{il}{h}-\frac{\varphi}{h}+i\frac{\psi_1}{h}}
b(\overline{a_1+r_1})dx \\
=\int_{\partial\Omega_+,\epsilon_0}
(\partial_\nu u)e^{-\frac{1}{h}(\varphi-i\psi_1)}\overline{(a_1+r_1)} dS.
\end{array}
\end{equation}
The second term of the LHS is what is different from (\ref{use.10}). Because
of (\ref{co.6}) this term goes to $0$ as $h$ goes to zero, since
\[
|e^{\frac{il}{h}-\frac{\varphi}{h}+\frac{i\psi_1}{h}}|=e^{-\frac{k(x)}{h}},
\]
and $q,b,a_1$, are bounded and $\|r_1\|_{H^0}\to 0$, $h\to 0$. Therefore
we get, instead of (\ref{use.14}), 
\[
\begin{array}{c}
\left|\int_{\partial\Omega_+,\epsilon_0} (\partial_\nu
u)e^{-\frac{1}{h}(\varphi-i\psi_1)} (\overline{a_1+r_1})
S(dx)\right|^2 \\
\le \frac{Coh}{\epsilon} \|a_1+r_1\|^2_{\partial\Omega_+,\epsilon_0}
\|e^{-\frac{\varphi}{h}}qu_2\|^2.
\end{array}
\]
The previous estimates imply that 
\[
\|e^{-\frac{\varphi}{h}}q u_2\|,\
\|a_1+r_1\|_{\partial\Omega_+,\epsilon_0}={\cal O}(1).
\]
Consequently the RHS of (\ref{co.12}) tends to $0$ as $h\to 0$ and we get
(\ref{use.16}) as before, namely
\begin{equation}\label{co.13}
\int_\Omega q(x)a_2(x)\overline{a_1(x)} e^{if(x)}dx=0.
\end{equation}

Now the arguments of Section \ref{SectionPf} imply that
$q=0$ finishing the proof of Theorem \ref{Th0.1}.
\hfill$\Box$
\end{document}